\documentclass[10pt,a4paper,twoside,reqno]{elsart}

\usepackage[T1]{fontenc}
\usepackage[latin1]{inputenc}
\usepackage[english,francais]{babel}

\usepackage{lscape}     
\usepackage[french]{varioref} 

\usepackage{ifthen}      

\usepackage{enumerate}   

\usepackage{array}       
\usepackage{multirow}    
\usepackage{booktabs}

\usepackage{graphicx}    
\usepackage{rotating}   
\usepackage{color}
\usepackage{epsfig}
\usepackage{subfig}

\usepackage{float}

\usepackage{amsmath}     
\usepackage{amssymb}     
\usepackage{mathrsfs}    

\floatstyle{ruled}
\newfloat{algorithme}{hbt}{loalg}[section]
\floatname{algorithme}{\textbf{Algorithme}}

\newcommand{\abs}[1]{\ensuremath{\left|#1\right|}}
                                
\newcommand{\R}{\ensuremath{\mathbb{R}}}

\newcommand{\dr}{\ensuremath{\partial}}
\newcommand{\dd}{\ensuremath{\mathrm{d}}}

\newcommand{\indi}[1]{\ensuremath{1\!\!1_{#1}}}
\renewcommand{\e}[1]{\ensuremath{\mathrm{e}^{#1}}}
\newcommand{\D}{\ensuremath{\mathscr{D}}}

\newcommand{\BL}{\ensuremath{\mathrm{cl}}}
\newcommand{\id}{\ensuremath{\mathrm{I}}}

\newcommand{\inc}{\ensuremath{\mathrm{inc}}}

\newcommand{\fonction}[5][]{\ifthenelse{\equal{#1}{}}%
  {\ensuremath{%
    \begin{array}[t]{@{}cll}
      #2&\longrightarrow&#3\\
      #4&\longmapsto&#5
    \end{array}}
  }
  {\ensuremath{%
    \begin{array}[t]{@{}lcll}
      #1\ :&#2&\longrightarrow&#3\\
      &#4&\longmapsto&#5
    \end{array}}
  }
}

\AtBeginDocument{}

\renewcommand{\phi}{\varphi}
\renewcommand{\epsilon}{\varepsilon}
\renewcommand{\S}{\ensuremath{\mathcal{S}}}

\renewcommand{\O}[2][]{\ifthenelse{\equal{#1}{}}%
   {\ensuremath{\mathcal{O}\left(#2\right)}
   }
   {\ensuremath{\underset{#1}{\mathcal{O}}\left(#2\right)}
   }
}

\renewcommand{\o}[2][]{\ifthenelse{\equal{#1}{}}%
   {\ensuremath{\text{\scriptsize $\mathcal{O}$}\left(#2\right)}
   }
   {\ensuremath{\underset{#1}{\text{\scriptsize $\mathcal{O}$}}\left(#2\right)}
   }
}

\renewcommand{\i}{\ensuremath{\mathrm{i}}}

 {\everymath{\displaystyle\everymath{}}\array}%
 {\endarray}

\journal{the Acad\'emie des sciences}
\begin{document}
\centerline{}
\begin{frontmatter}




%
\selectlanguage{francais}
\title{Mod\`ele \'electromagn\'etique d'objet dissimul\'e}



\author[authorlabel1]{Jean-Baptiste Bellet},
\ead{bellet@cmap.polytechnique.fr}
\author[authorlabel2]{G\'erard Berginc}
\ead{gerard.berginc@fr.thalesgroup.com}

 \address[authorlabel1]{CMAP, \'Ecole Polytechnique, Route de Saclay, 91128 Palaiseau Cedex}
 \address[authorlabel2]{Thales Optronique, 2, Avenue Gay Lussac CS 90502, 78995 \'Elancourt Cedex}



 \medskip
 \selectlanguage{francais}
 \begin{center}
 {\small Re\c{c}u le *****~; accept\'e apr\`es r\'evision le +++++\\
 Pr\'esent\'e par £££££}
 \end{center}

\begin{abstract}
\selectlanguage{francais}
Nous \'elaborons un mod\`ele de propagation d'ondes \'electromagn\'etiques dans un milieu, inhomog\`ene, avec une couche rugueuse, et dissimulant un objet. Nous obtenons un milieu effectif, puis  nous r\'esolvons le probl\`eme par \'equations int\'egrales.
{\it Pour citer cet article~: J.-B. Bellet, G. Berginc, C. R.
Acad. Sci. Paris, Ser. * *** (2010).}
\vskip 0.5\baselineskip

\selectlanguage{english}
\noindent{\bf Abstract}
\vskip 0.5\baselineskip
\noindent
{\bf Electromagnetic model of a hidden object. }
We model propagation of electromagnetic waves in a medium, which is inhomogeneous with a rough layer, and which hides an object. We first get an effective medium, and then we solve the problem by integral equations.
{\it To cite this article: J.-B. Bellet, G. Berginc, C. R.
Acad. Sci. Paris, Ser. * *** (2010).}
\end{abstract}


\end{frontmatter}

\selectlanguage{francais}

\section*{Introduction}
Nous sommes motiv\'es par l'imagerie laser d'un objet dans un milieu fluctuant. Ceci est adapt\'e pour la d\'etection d'une tumeur enfouie dans le derme, ou d'un engin explosif improvis\'e camoufl\'e. Nous r\'esolvons un probl\`eme mod\`ele de propagation des ondes \'electromagn\'etiques dans un milieu contenant un demi-espace avec des inhomog\'en\'eit\'es, surmont\'e d'une couche rugueuse. La couche peut repr\'esenter la couche corn\'ee de la peau, ou un filet de camouflage. Le demi-espace peut constituer l'int\'erieur de la peau, ou un milieu de camouflage. Nous nous pla\c{c}ons dans un contexte bi-dimensionnel p\'eriodique. Enfin, nous consid\'erons un cas asymptotique multi-\'echelles : la couche est mince par rapport \`a la longeur d'onde, et les inhomog\'en\'eit\'es sont petites par rapport \`a la couche mince.
\\
Ainsi, on mod\'elise la couche sup\'erieure par la couche mince p\'eriodique $\D^\BL_\xi=\{(x_1,x_2):0<x_2<\xi f(x_1/\xi)\}$, o\`u $f>0$ est une fonction 1-p\'eriodique, et $\xi>0$ est un petit param\`etre (par rapport \`a la longueur d'onde). Cette couche mince est d\'elimit\'ee par les interfaces $\gamma_0=\{x_2=0\}$ et $\gamma_\xi=\{(x_1,\xi f(x_1/\xi)),x_1\in\R\}$. La sous-couche est le demi-espace $\D^-=\{x_2<0\}$. Le milieu ext\'erieur est $\D^+_\xi=\{x_2>\xi f(x_1/\xi)\}$. Le milieu sup\'erieur $\D^+_\xi$ et la couche $\D^\BL_\xi$ sont homog\`enes. Le milieu inf\'erieur $\D^-$ est un milieu $\xi Y$-p\'eriodique contenant des petites inclusions, o\`u la cellule de r\'ef\'erence $Y=(0,\ell_1)\times(0,\ell_2)$ ($\ell_1,\ell_2>0$) contient une inhomog\'en\'eit\'e $B$. Bien que la taille caract\'eristique de la p\'eriode soit encore not\'ee $\xi$ dans $\D^-$, nous la supposons \^etre \`a une \'echelle inf\'erieure \`a la taille de la couche.  Dans la cellule $Y$, la perm\'eabilit\'e et la permittivit\'e sont respectivement : $\mu_Y(y)=\mu_B\indi{B}(y)+\mu\indi{Y\setminus\overline{B}}(y)$, $\epsilon_Y(y)=\epsilon_B\indi{B}(y)+\epsilon\indi{Y\setminus\overline{B}}(y)$, o\`u $\mu_B,\mu,\epsilon_B,\epsilon>0$ sont des constantes. Dans tout le domaine $\R^2$, la perm\'eabilit\'e et la permittivit\'e sont respectivement : $\mu_\xi(x)=\mu^+\indi{\D^+_\xi}(x)+\mu^\BL\indi{\D^\BL_\xi}(x)+\mu_Y(x/\xi)\indi{\D^-}(x),$ $\epsilon_\xi(x)=\epsilon^+\indi{\D^+_\xi}(x)+\epsilon^\BL\indi{\D^\BL_\xi}(x)+\epsilon_Y(x/\xi)\indi{\D^-}(x)$, o\`u $\mu^+,\epsilon^+>0$ sont des constantes et o\`u l'on a prolong\'e par p\'eriodicit\'e les coefficients dans $\D^-$, et $\indi{\cdot}$ est la fonction indicatrice. On illumine dans la partie sup\'erieure du domaine par un laser, \emph{i.e.} une onde plane $u_\inc(x)=\e{\i k^+ \hat{\theta}\cdot x}$. Cette onde est de pulsation $\omega>0$, de nombre d'onde $k^+=\omega\sqrt{\epsilon^+\mu^+}$, d'angle d'incidence $\hat{\theta}=(\hat{\theta}_1,\hat{\theta}_2)$ avec $\hat{\theta}_2<0$ et $\abs{\hat{\theta}}=1$. La longueur d'onde est $\lambda^+:=2\pi/k^+$, suppos\'ee grande par rapport \`a la taille de la couche mince $\xi$. Le probl\`eme mod\`ele consiste alors en l'\'etude de la propagation de l'onde r\'esultante $u_\xi$ dans le milieu :
\begin{equation}
\nabla\cdot \frac{1}{\mu_\xi}\nabla u_\xi
+\omega^2 \epsilon_\xi u_\xi=0,\hfill\text{ dans }\R^2,\\
\label{eq:pb_modele}
\end{equation}
avec en plus une condition de radiation. Cette derni\`ere s\'electionne la solution du physique du probl\`eme ; en particulier, elle veille \`a ce que l'onde r\'efl\'echie $u_\xi-u_\inc$ dans le milieu sup\'erieur $\D^+_\xi$  et l'onde transmise $u_\xi$ dans le milieu inf\'erieur $\D^-$ soient des ondes sortantes.
Enfin, une tumeur dans le derme, ou une bombe camoufl\'ee, est mod\'elis\'ee par une inclusion $D\subset \D^-$,  de permittivit\'e $\epsilon_D>0$ et de perm\'eablit\'e $\mu_D>0$, enfouie dans la partie inf\'erieure $\D^-$. On supposera que $D$ est suffisamment en profondeur par rapport \`a la taille de la couche, et que $D$ est suffisamment grande devant les inhomog\'en\'eit\'es du milieu inf\'erieur.
\\
Notre probl\`eme mod\`ele, certes simpliste, pr\'esente l'avantage d'\^etre purement \'electromagn\'etique, et de prendre en compte la rugosit\'e de l'interface. Il est nouveau ; pour le r\'esoudre, nous proposons une d\'emarche originale, en adaptant pas \`a pas des techniques classiques. Son caract\`ere multi-\'echelles et p\'eriodique conduit \`a l'approcher par un mod\`ele effectif qui lui est \'equivalent, par des techniques asymptotiques usuelles \cite{ALL.CON,ABB.DIFF,CIU.APP}. Par hypoth\`eses sur l'objet $D$, son ajout  dans le milieu p\'eriodique est \'equivalent \`a son ajout dans le milieu effectif. La solution du probl\`eme initial en pr\'esence de $D$ est alors calcul\'ee par la m\'ethode des \'equations int\'egrales \cite{NED.ACO,AMM.POL}, dans le milieu effectif.


\section{Milieu effectif}

Le caract\`ere multi-\'echelles p\'eriodique du probl\`eme \eqref{eq:pb_modele} sugg\`ere d'appliquer les m\'ethodes asymptotiques usuelles. La diff\'erence d'ordre des \'echelles permet d'effectuer l'\'etude en deux temps. Tout d'abord, de fa\c{c}on similaire \`a un calcul men\'e dans \cite{ALL.CON}, on homog\'en\'eise le milieu inf\'erieur $\D^-$ en tronquant un d\'eveloppement \`a deux \'echelles. Puis, de fa\c{c}on analogue \`a \cite{ABB.DIFF,CIU.APP}, on \'elabore des conditions de transmission, de type imp\'edance g\'en\'eralis\'ee, permettant d'approcher l'effet de la couche mince. Cel\`a s'obtient en \'ecrivant un d\'eveloppement asymptotique \`a l'aide de correcteurs de couche limite d\'ecroissant exponentiellement. On trouve alors que la solution du probl\`eme mod\`ele~\eqref{eq:pb_modele} peut \^etre approch\'ee par $U$ solution d'un probl\`eme dans deux demi-espaces homog\`enes, avec des conditions de transmission \`a l'interface :
\begin{equation}
\left\{
\begin{array}{l}
\frac{1}{\mu^+}\Delta U
+\omega^2 \epsilon^+ U=0,\text{ dans }\D^+,\hfill\nabla\cdot \mathcal{A}\nabla U+\omega^2 \epsilon^- U=0,\text{ dans }\D^-,\\
\left[U\right]=\psi\cdot\left.\nabla U\right|_{\gamma_0^+}, \hfill\text{ sur }\gamma_0,\\
\frac{1}{\mu^+}\left.\dr_{x_2}U\right|_{\gamma_0^+}-(0,1) \mathcal{A}\left.\nabla U\right|_{\gamma_0^-}=\phi_1\cdot\left.\dr_{x_1}(\left.\nabla U\right|_{\gamma_0^+})\right|_{\gamma_0}\\
\hfill+\phi_2\cdot \left(
                    \dr^2_{x_2 x_1} U|_{\gamma_0^+},
                     (-{k^+}^2-\dr^2_{{x_1}^2}) U|_{\gamma_0^+}
                   \right)^T
+\phi_3\left.U\right|_{\gamma_0^+}, \quad\text{ sur }\gamma_0.
\end{array}\right.
\label{eq:pb_effectif}
\end{equation}
Les param\`etres effectifs du milieu inf\'erieur sont issus de l'homog\'en\'eisation, donn\'es par $\mathcal{A}=\frac{1}{\abs{Y}}\int_Y \frac{1}{\mu_Y(y)}\left(\id-\nabla\chi(y)\right)\dd y$, et $\epsilon^-=\frac{\abs{B}}{\abs{Y}}\epsilon_B+(1-\frac{\abs{B}}{\abs{Y}})\epsilon$, o\`u $\chi$ est la solution $Y$-p\'eriodique \`a moyenne nulle du probl\`eme de cellule $\nabla_y\cdot\frac{1}{\mu_Y(y)}\nabla_y \chi(y)=\nabla_y\cdot \frac{1}{\mu_Y(y)}\id, y\in Y$.
Pour les conditions de transmission, les coefficients sont donn\'es par :
$s=\mathcal{A}_{21}$, $\psi=\xi\psi_0$, $\phi_2=\xi \left(\frac{1}{\mu^\BL}-\frac{1}{\mu^+}\right)\int_0^1 f(0,1)^T$, $\phi_3=\xi \omega^2\left(\epsilon^\BL-\epsilon^+\frac{\mu^+}{\mu^\BL}\right)\int_0^1 f$, $\phi_1/\xi=-s\int_{\Gamma_0}\left.\Psi_0\right|_{\Gamma_0}\dd\sigma+\left(\frac{1}{\mu^+}-\frac{1}{\mu^\BL}\right)\int_{\Gamma_1}\nu_1 \left.\Psi_0\right|_{\Gamma_1} \dd\sigma$, et  $(\Psi_0,\psi_0)$ est l'unique solution du probl\`eme de cellule :
\begin{equation*}
\left\{\begin{array}{l}
\nabla\cdot \tilde{\mathcal{A}}\nabla \Psi_0=\left(\frac{1}{\mu^\BL}-\frac{1}{\mu^+}\right)\nu \delta_{\Gamma_1}-\left(\frac{1}{\mu^+}-\frac{1}{\mu^\BL}\right)\nu \delta_{\Gamma_0},\hfill \text{ dans }(0,1)\times\R, \\
y_1\longmapsto \Psi_0(y_1,\cdot),\hfill 1\text{-p\'eriodique},\\
\Psi_0\longrightarrow 0,  y_2\to-\infty;\hfill \Psi_0\longrightarrow \psi_0, y_2\to+\infty,
\end{array}\right.
\end{equation*}
o\`u par remise \`a l'\'echelle $y:=x/\xi$, $\Gamma_1$ et $\Gamma_0$ sont les images respectives de $\gamma_\xi$ et $\gamma_0$, et $\tilde{\mathcal{A}}(y)=\mathcal{A}\indi{\D^-}(x)+\frac{1}{\mu^\BL}\indi{\D^\BL_\xi}(x)+\frac{1}{\mu_+}\indi{\D^+_\xi}(x)$.\\
La solution $U$ du probl\`eme effectif~\eqref{eq:pb_effectif}, et la fonction de Green $G$ associ\'ee, peuvent \^etre calcul\'ees par analyse de Fourier tangentielle \`a l'interface $\gamma_0$. Pour $G$, ceci n\'ecessite alors l'\'evaluation d'int\'egrales de Sommerfeld.



\section{Probl\`eme mod\`ele d'objet dissimul\'e}
Nous avons transform\'e le probl\`eme mod\`ele~\eqref{eq:pb_modele} en le probl\`eme effectif~\eqref{eq:pb_effectif}. Ajoutons une inclusion $D\subset \D^-$ dans la partie inf\'erieure du milieu, de permittivit\'e $\epsilon_D>0$ et de perm\'eablit\'e $\mu_D>0$.  On suppose que $D$ est suffisamment en profondeur par rapport \`a la taille de la couche, et que $D$ est grande devant les inhomog\'en\'eit\'es du milieu inf\'erieur. Sous ces hypoth\`eses, l'analyse asymptotique permet,  formellement, de remplacer le milieu de fond par le milieu effectif pr\'ec\'edent. Pour r\'esoudre ce probl\`eme, nous utilisons la m\'ethode des \'equations int\'egrales. On note $G_D$ la fonction de Green de l'anomalie ; on d\'efinit  les potentiels de simple couche par $\S\phi(x):=\int_{\dr D} G(x,y)\phi(y)\dd\sigma(y)$, $\S_D\phi(x):=\int_{\dr D} G_D(x,y)\phi(y)\dd\sigma(y)$. L'ajout de $D$ transforme $U$ en $u$ par convolution avec les noyaux de Green \cite{NED.ACO,AMM.REC,AMM.POL} :
\begin{equation*}
  u=U+\S\phi\text{ dans }\overline{D}^c, u=\S_D\psi\text{ dans }D,
\end{equation*}
o\`u  $(\phi,\psi)\in L^2(\dr D)\times L^2(\dr D) $ est solution d'un syst\`eme d'\'equations int\'egrales permettant d'assurer les conditions de transmission sur $\dr D$. Ce dernier est r\'esolu num\'eriquement par la m\'ethode des \'el\'ements finis de fronti\`ere.



\section*{Remerciements}
Les auteurs remercient tr\`es chaleureusement Habib Ammari pour ses conseils et suggestions.

\bibliographystyle{elsart-num-sort}
\bibliography{biblio}
\addcontentsline{toc}{section}{R\'ef\'erences}



\end{document}